\documentclass[12pt]{article} 
\usepackage{amssymb,amsmath} 
\newtheorem{theorem}{Theorem} 
\newtheorem{lemma}{Lemma} 
\newtheorem{proposition}{Proposition} 
\newtheorem{corollary}{Corollary}

\begin{document} 
\title{Legendrian Submanifold Path Geometry}
\author{ Sung Ho Wang \\
Department of Mathematics\\
Duke University\\
e-mail : ship@math.duke.edu\\ } 
\date{April 15, 2000} 

\maketitle

\section{Introduction}

In [Ch1], Chern gives a generalization of projective geometry
by considering foliations on the Grassman bundle of $p$-planes
$Gr(p,R^{n}) \to R^{n}$ by $p$-dimensional 
submanifolds that are integrals of the canonical
contact differential system. The equivalence method yields
an $\mathfrak{sl}(n+1,R)$-valued Cartan connection whose curvature captures the 
geometry of such foliation. In the flat case, the space of
 leaves of the foliation is a second order homogeneous space
 [Br2].

[Ch2] deals with the geometry of the foliation of $Z^4 \to Y^3$, where $Z$
is the bundle of Legendrian line elements over a contact 
threefold $Y$, by canonical lifts of Legendrian curves, or 
equivalently, the geometry of
3-parameter families of curves in the plane. An $\mathfrak{sp}(2,R)$ -valued
Cartan connection plays the role of projective connection.

A generalization of [Ch2] to 4-parameter family of 
curves in the plane leads to a geometric realization of 
some exotic holonomies in dimension four [Br1].

In this paper, we generalize [Ch2] to higher dimensions. 
Let $Z \to Y^{2n+1}$ be the bundle of Legendrian $n$-planes
over a contact manifold $Y$. We consider a foliation  
of $Z$ by  canonical lifts of Legendrian submanifolds,
which we call a \emph{Legendrian submanifold path geometry}. Note
that a path
in this case is a Legendrian $n$-fold. The 
equivalence method provides an $\mathfrak{sp}(n+1, R)$-valued Cartan 
connection form that captures the geometry of such 
foliation. In the flat case, the space $X$ of leaves
of the foliation is again a second order homogeneous 
space. The  prolonged structure equation of
this second order homogeneous space is in turn that 
of $Sp(n+1, R)$, which explains the appearance of a $\mathfrak{sp}(n+1, R)$ 
-valued Cartan connection form. In fact, we may consider a contact
manifold $Y$ endowed with a Legendrian submanifold path geometry structure as a
 union of infinitesimal homogeneous spaces 
\[ Sp(n+1, R) \to RP^{2n+1} \]
connected by the Cartan connection mentioned above. As a by 
product, this gives
a geometric realization of the Lie algebra $\mathfrak{sp}(n+1, R)$
 as the 
symmetry vector fields of a family of Legendrian submanifolds
in a $(2n+1)$ -dimensional contact manifold. 

After a short analysis of the structure equations associated
to a general Legendrian submanifold path geometry, 
two special cases are considered. The first case is 
characterized by having a well defined conformal class of
symmetric $(n+1)$ differentials on the space of leaves of
the foliation $X$. The vanishing
of this symmetric differential represents a necessary 
condition for the contact of neighboring Legendrian leaves.
A double fibration naturally arises, and we give a dual
interpretation of the contact manifold $Y$ in terms of 
$X$. The $G$-structure induced on $X$ gives an example
of a classical non-metric irreducible
holonomy $GL(n+1,R)$ with representation on $sym^2(R^{n+1})$.
 
It is well known that the normal projective connection
uniquely associated to a torsion free affine connection
captures the geometry of paths defined by the geodesics
of the affine connection.  The second case is a direct
generalization of this to Legendrian submanifold path geometry. We consider a 
\emph{Legendrian} connection on the contact hyperplane
vector bundle over $Y$ whose \emph{geodesic} Legendrian
submanifolds give rise to a Legendrian submanifold path geometry. 
There exists a unique \emph{normal symplectic} connection
associated to a Legendrian connection such
that any other Legendrian connection with the 
equivalent Legendrian submanifold path geometry is a 
section of the normal symplectic connection. An analysis
of the normal symplectic connection shows in fact the family
of geodesic isotropic $k$-folds for $1 \leq k \leq n-1$
 with respect to a Legendrian connection is 
also an invariant of the normal symplectic connection. 

$RP^{2n+1}$, as a quotient space of $Sp(n+1, R)$, carries a 
Legendrian submanifold path geometry, which is flat. For a nonflat
example with symmetry, consider a hypersurface $M^n$ in the
$(n+1)$ -dimensional space form $\bar{M}_c^{n+1}$, $c=1, 0,$
or $-1$, without any extrinsic symmetry. 
The images of $M$ under the motion by Iso($\bar{M}_c^{n+1}$), when lifted to
$Gr(n, \bar{M}_c^{n+1} )$, generates an
$N=\frac{1}{2}(n+1)(n+2)$-parameter family of
Legendrian submanifolds. Since  Iso($\bar{M}_c^{n+1}$) does not
arise as a subgroup of $Sp(n+1,R)$, it is not equivalent 
to the flat example.  

Similar constructions are likely to work for other
(irreducible) second order homogeneous spaces [Br2]. For instance,
in the holomorphic category for simplicity, a manifold
with $CO(V)$ structure has, as its dual, a manifold with
$GL(W)$ structure with representation $\bigwedge^2(W)$ where
$V$ and $W$ are vector spaces of suitable dimensions. 
The corresponding Cartan connection form would be 
$\mathfrak{so}(m,\mathbb{C})$-valued for suitable $m$.
Two exceptional cases, $\mathbb{C}^* Spin(10, 
\mathbb{C})$ on $\mathbb{S}_+$ and $\mathbb{C}^* 
E_6^{\mathbb{C}}$ on $\mathbb{C}^{27}$, would yield 
Cartan connection forms 
with values in $\mathfrak{e}_6^{\mathbb{C}}$ and 
$\mathfrak{e}_7^{\mathbb{C}}$
, respectively, for the associated geometries.

All the arguments remain valid when we replace real and
smooth by complex and holomorphic. In fact, a \emph{real}
Legendrian submanifold path geometry can be considered as a split real
form of a \emph{complex} one. In analogy with the geometry
of real hypersurfaces in $\mathbb{C}^{n}$ considered as a real form 
of a complex hypersurface path geometry via Segre families [Fa] [ChM],
it would be interesting to study other possible real 
forms of complex Legendrian submanifold path geometry.

We shall agree that all the Latin indices $i, j, k$ run from
$1$ to $n$, and, for simplicity, that $n \geq 2$.

We would like to thank Prof. Robert Bryant for his guidance and 
support throughout this work.

\newpage

\section{Legendrian Submanifold Path Geometry}

\subsection{Definition}

Let $Y$ be a $(2n+1)$ -dimensional manifold with a contact
structure, i.e., a differential 1-form $\theta_0$ defined
up to multiplication by a nonzero function with the property 
\begin{equation}
 \theta_0 \wedge ( d\theta_0 )^n \ne 0. 
\end{equation}

The differential system on $Y$
locally generated by $\theta_0$ and $d\theta_0$
is called a \emph{contact differential system.}
It is known that at each point of $Y$, a subspace of 
the tangent space of $Y$ that is integral to the contact
differential system is of 
dimension at most $n$ and such a subspace is called a Legendrian
$n$-plane. An integral $n$-dimensional submanifold
is similarly called a Legendrian $n$-fold or Legendrian 
submanifold.  Let $Z \to Y$ be the associated
bundle of Legendrian $n$-planes, which we may regard as
the first prolongation of the contact differential system
on $Y$. The theorem of Darboux on the local normal form for the 
contact structures provides local coordinates 
 $ \{\; x^i, u, p_i, p_{ij} \; \}$
 on $Z$ such that $\{ \; x^i, u, p_i \; \} $ is a local 
coordinate system on $Y$ and so that the first prolongation 
of the contact system to $Z$ is generated by
\begin{align} 
\theta_0 &=  du -  \sum_{k=1}^n p_k \, dx^k  \notag \\ 
\theta_i &=  dp_i -\sum_{k=1}^n p_{ik} \, dx^k, 
                        \; i = 1, .. \, n. \notag 
\end{align}
This means that an $n$-dimensional integral submanifold
of the differential system generated by
\[ \{ \; \theta_0,  \theta_1, .. , \theta_n, 
             d\theta_1, .. , d\theta_n \; \} \]
on which 
$ dx^1 \wedge ..  \wedge dx^n \ne 0 $
is the canonical lift of a Legendrian $n$-fold in $Y$.
In fact, given any point in a Legendrian submanifold 
in $Y$, there exists a local coordinate system
$\{\; x^i, u, p_i \;\}$ of $Y$, in the neighborhood of the given
point, such that the given
Legendrian submanifold is defined by the equations
  $\{ \, u=0, \, p^1=0, \; .. \; , \, p^n=0 \, \}$,
the contact form
$\theta_0$ is a multiple of
\[ du -  \sum_{k=1}^n p_k \, dx^k, \] 
and $dx^1 \wedge ..  \wedge dx^n \ne 0 $ on the given
Legendrian submanifold, at least in a neighborhood of 
the given point.

The problem we are interested in is a local geometry of
the foliation of $Z$, necessarily transversal to the fibers
of the projection $Z \to Y$, by $n$-dimensional 
submanifolds that are lifts of Legendrian 
$n$-folds in $Y$. In terms of $Y$, this means to each 
Legendrian $n$-plane (possibly only those in an open set of the set of Legendrian
$n$-planes) at a point, there exists a unique
Legendrian submanifold of this family tangent to the 
given Legedrian $n$-plane. Let $X$ denote the space
of the leaves of the foliation, which we always assume
to be a nice $N=\frac{1}{2}(n+2)(n+1)$ -dimensional 
manifold. 

We could also describe this geometry as
the geometry of a nondegenerate $N$-parameter family of 
Legendrian submanifolds in $Y$. Here, an $N$-parameter family
 of Legendrian submanifolds is said to be non degenerate if the
following is true: Let X$^N$ be the parameter space
and $l$ : X$\times R^n \to Y$ the $N$-parameter family
of Legendrian immersions of $R^n$ to $Y$. We require that the
associated lift $\hat{l}$ : X$\times R^n \to Z$ be
a (local) diffeomorphism. Roughly, this means that the union of the lifts 
of the $N$-parameter family of Legendrian submanifolds 
to $Z$ fills out an open set in $Z$. 

\vspace{1pc}
\noindent \textbf{Definition} 
\textnormal{Legendrian Submanifold Path Geometry} 

Let $Z \to Y^{2n+1}$ be the bundle of Legendrian $n$-planes
of a contact manifold $Y$. \emph{Legendrian submanifold path geometry} is the
geometry of  foliations on $Z$, up to the diffeomorphism
of $Z$ induced from the contact transformation of $Y$,
whose leaves are the cannonical lifts
of Legendrian submanifolds in $Y$. Locally, this is 
equivalent to the geometry of a non degenerate 
$N=\frac{1}{2}(n+1)(n+2)$-parameter family of Legendrian 
submanifolds in $Y$ up to contact transformations.
\vspace{1pc}

\subsection{Structure Equations}

In this section, we give a brief analysis of the structure
equation for general Legendrian submanifold path geometry.

Let $\mathcal{I}$ be the Frobenius system on $Z$
describing the given Legendrian submanifold path geometry. 
The local normal form theorem above shows that
$\mathcal{I}$ is locally generated by the
following 1-forms.
\begin{align} 
\theta_0 &=  du - \sum_{k} p_k \, dx^k   \\ 
\theta_i &=  dp_i -\sum_k p_{ik} \, dx^k, 
                        \; i = 1, .. \, n \notag  \\ 
\Theta_{ij} &= dp_{ij} - \sum_k F_{ijk} \, dx^k,
                        \; i , j = 1, .. \, n, \notag 
\end{align} 
where $\{\; x^i, u, p_i, p_{ij} \; \} $ is the local coordinate
system mentioned above and $\{ F_{ijk} \}$ is a set of
functions locally defined on $Z$.
These differential forms form a subset of a coframe 
$\{ \, \omega^i = dx^i, \theta_0, \theta_i, \Theta_{ij} \}$
of $Z$, defined up to the diffeomorphisms of $Z$ induced 
from the contact transformations of $Y$, with
\begin{equation}
\mathcal{I} = \{ \, \theta_0, \theta_i, \Theta_{ij} \, \}
\end{equation}
and
\begin{align}  
 d\theta_0 &\equiv - \sum_k \theta_k \wedge \omega^k  
               \quad \mod \quad \theta_0   \\ 
 d\theta_i &\equiv - \sum_k \Theta_{ik} \wedge \omega^k
     \; \; \; \mod \quad \theta_0, \theta   \notag \\
 d\Theta_{ij} &\equiv 
    \quad 0 \quad \quad \mod \quad \theta_0, \theta, \Theta,  \notag 
\end{align}
where$\mod \; \theta \,$ means$\mod \; \theta_1, .. \, 
\theta_n \,$, and similarly for$\mod \; \Theta \,$. 

Given such a differential system $\mathcal{I}$ on $Z$, the 
set of all coframes $\{ \, \omega^i, \theta_0, \theta_i, 
\Theta_{ij} \}$
on $Z$ satisfying (3) and (4) form a 
$H_1 \subset Gl(N+n,R)$ bundle $F_1$ over $Z$, where $H_1$ is 
the subgroup whose induced action on $T^* Z$ preserves (3) and (4).
Equivalently, the principal right $Gl(N+n,R)$ coframe bundle 
can be reduced to a $H_1$ subbundle via $\mathcal{I}$.  
It can be easily shown the right action of $H_1$ on the 
tautological forms of $F_1$, which are by definition the 
restriction to $F_1$ of the $Gl(N+n,R)$ equivariant 
$R^{N+n}$-valued tautological 1-form on the 
principal $Gl(N+n,R)$ bundle, is as follows.
\begin{align}  
  \theta_0^* &= \lambda \theta_0 \\ 
  \theta^*   &= A( \theta + b \, \theta_0)  \notag \\
  \Theta^*   &= \frac{1}{\lambda} A (\Theta 
           + e \theta_0 + \sum_i A_i \, \theta_i) A^{t}
               \notag \\ 
  \omega^*   &\equiv \lambda A^{-t} ( \omega + c \, \theta )
              \quad \mod \theta_0,  \notag             
\end{align} 
where we denote 
\begin{align}
\theta^t &= ( \, \theta_1, .. \, , \theta_n \, ) \notag \\
\Theta &= ( \, \Theta_{ij} \, ) = 
          ( \, \Theta_{ji} \, ) \notag \\
\omega^t &= ( \, \omega^1, .. \, , \omega^n \, ), \notag 
\end{align}
and $\lambda \ne 0$, $A \in Gl(n,R)$, 
$A_i, e, c \in sym^2(n,R)$, and $b \in R^n$. Here $sym^2(n,R)$
denotes the set of $n$ by $n$ symmetric matrices. We used 
the same notation $\omega^i, \theta_0, \theta_i, 
\Theta_{ij}$ as above to denote the corresponding 
tautological forms.

Thus on $F_1$, we have the following structure equations.
\begin{align}
d\theta_0 &= -\rho \wedge \theta_0 - 
             \theta^t \wedge \omega + T_{00} 
                    \wedge \theta_0  \notag  \\
d\theta &= -\beta \wedge \theta_0 -\alpha \wedge \theta 
             - \Theta \wedge \omega + \sum_k T_{11}^k
    \wedge \theta_k + T_{10} \wedge \theta_0 \notag \\
d\Theta &= - \epsilon \wedge \theta_0
           - \sum_k \pi^k \wedge \theta_k 
           -( \alpha - \frac{\rho}{2} ) \wedge \Theta 
           + \Theta \wedge ( \alpha - \frac{\rho}{2} )^t
                    \notag \\
        & \quad + \sum_k T_{21}^k \wedge \theta_k + 
           T_{20} \wedge \theta_0 + T_0 \notag \\
d\omega &= - \mu \wedge \theta_0 - \gamma \wedge \theta 
           +( \alpha^t - \rho ) \wedge \omega
           + T^1_0 \wedge \theta_0
           + T^{11}. \notag
\end{align}
Here $\rho$ is a scalar 1-form, $\alpha$ is a  $gl(n,R)$
-valued 1-form, $\pi^k, \epsilon, \gamma$ are symmetric  
$gl(n,R)$ -valued 1-forms, and $\beta, \mu$ are
$R^n$ (column) valued 1-forms. These  
are the  pseudo connection forms on $F_1 \to Z$, and 
$T_{00}$, $T_{11}^k$, $T_{10}$, $T_{21}^k$, $T_{20}$, $T_0$, 
 $T^1_0$, and $T^{11}$ represent the torsion 
of this pseudo connection. The pseudo connection forms are
not uniquely defined. By modifying the pseudo connection 
forms, we may reduce the torsion to the following 
simple form.

\begin{proposition}
There exists a pseudo connection on $F_1$ for which the
structure equation takes the following form.
\begin{align}
d\theta_0 &= -\rho \wedge \theta_0 - 
             \theta^t \wedge \omega   \\
d\theta &= - \beta \wedge \theta_0 -\alpha \wedge \theta 
             - \Theta \wedge \omega \notag \\
d\Theta &= - \epsilon \wedge \theta_0 
           - \sum_k \pi^k \wedge \theta_k
           -( \alpha - \frac{\rho}{2} ) \wedge \Theta 
           + \Theta \wedge ( \alpha - \frac{\rho}{2} )^t
           + T_0 \notag \\
d\omega &=  - \mu \wedge \theta_0 - \gamma \wedge \theta 
            +( \alpha^t - \rho ) \wedge \omega,
            \notag
\end{align}
where $T_0 = \sum_{ij} \Theta_{ij} \wedge \tau^{ij}$ with each 
$\tau^{ij} = \tau^{ji}$ being symmetric $gl(n,R)$
 -valued 1-form satisfying $\tau^{ij} \wedge \omega = 0$.
\end{proposition}

\emph{Proof.} First, by modifying $\rho$, and $\beta$, we 
can absorb $T_{00}$, and $T_{10}$. Also, the second 
equation in  (4)
implies that we can arrange $T_{11}^k$'s to be $0$ by modifying
$\alpha$. Thus all the torsion terms in 
$d\theta_0$ and $d\theta$ can be absorbed, which we assume
from now on.

We modify $\mu$ to absorb $T^1_0$ and arrange $T^{11}$
to be quadratic in $\{ \, \omega, \theta, \Theta \, \}$.
Now, $d(d\theta_0) \equiv 0  \mod \theta_0$ gives
$ \theta^t \wedge T^{11} = 0$,
which means $T^{11}$ is of the form
\[ T^{11} = ( h_{ij} ) \wedge \theta \]
with $( h_{ij} )$ being a $gl(n,R)$ -valued 1 
form in  $\omega, \theta, \Theta$, which is 
not uniquely defined. Now, it is easily verified
 that by modifying the first row or column of 
$(\gamma_{ij}) = (\gamma_{ji}) $ 
and the representation $( h_{ij} )$, we may have $h_{1j}
 = 0$ and $h_{j1} = 0$ for $j = 1, .. \, n$ and 
$( h_{ij} ) = 0 \mod \omega, \theta_2, \theta_3, 
.. \, , \theta_n, \Theta$.
Hence by induction, we can absorb all of $T^{11}$
using $\gamma$.

Finally, we modify $\pi^k$, and $\epsilon$ to absorb
$T_{21}^k$, $T_{20}$ and arrange $T_0$ to be quadratic in
$\{ \, \omega, \Theta \, \}$. But $d(d\theta) \equiv 0
 \mod \theta_0, \theta$ gives
\[ T_0 \wedge \omega = 0, \]
which implies $T_0$ cannot have any quadratic terms in 
$\Theta$, and since $\mathcal{I}$ is Frobenius, it cannot
have any quadratic terms in $\omega$ either. $\square$

\vspace{1pc}
The torsion $T_0$, as it stands, is not an invariant of the
foliation. In fact the structure group $H_1$ acts on
 $T_0$. However, rather than 
continuing the analysis of equivalence problem directly, 
we examine a special case of a foliation motivated by [Ch2],
namely, that of quadric hypersurfaces in $R^{n+1}$.

\newpage
 
\section{Second Order Developables for 
Quadric Hypersurfaces in $R^{n+1}$ }
 
The local normal form (2) for $\mathcal{I}$ on $Z$ shows, at
least locally, we can identify $Z \to Y$ with 
$J^2(R^n, R^{n+1}) \to J^1(R^n, R^{n+1})$ 
and regard the geometry of the foliation as
the geometry of an $N$-parameter family of hypersurfaces in
$R^{n + 1}$ up to contact transformations. In this section,
 we take the simple example of $\mathcal{I}$
on $J^2(R^n, R^{n+1}) \cong Z$ defining the quadric 
hypersurfaces in $R^{n+1}$, 
\begin{align}
     u &= a_0 + \sum_i a_i x^i + \sum_{ij} \frac{1}{2} 
                    a_{ij} x^i x^j  \\
   p_i &= a_i + \sum_{j} a_{ij} x^j \notag \\
p_{ij} &= a_{ij} \notag
\end{align}
where $a_{ij}=a_{ji}$ and $\{u, x^i, p_i, p_{ij} \}$
is a local coordinate system of $Z$
introduced  earlier. Given the explicit
form of solutions, we may regard $\{\, a_0, a_i, a_{ij} \, \}$
as a local coordinate system on the space of solutions of 
$\mathcal{I}$. Also, it is easy to see that this family of
submanifolds in $J^1(R^n, R^{n+1}) \cong Y$ is 
non degenerate in the sense discussed earlier.

Consider a hypersurface S in $R^{n + 1}$
defined as the graph of a function $f$ in $n$ variables
\[ u = f( x_1, .. \, x_n ). \]
At each point of S, there exists a quadric hypersurface of
the form (7) that osculates the given hypersurface S up to
second order. The set of all such quadric hypersurfaces 
along S generically form an $n$-parameter family of solutions
to $\mathcal{I}$,
or an $n$ -dimensional submanifold $\mathcal{S}$ in the space 
of solutions. Coversely we may consider the 
original hypersurface S as the second order developable 
of the family $\mathcal{S}$.

From the construction, each quadric hypersurface 
in the family $\mathcal{S}$
has the point of contact with the given S,
$( \, u, x^1, .. \, x^n \, )$, at which
\begin{align}
 \theta_0 &= du - \sum_i p_{i} dx^i 
   = da_0 + \sum_i x^i da_i + \sum_{ij} \frac{1}{2} 
                   x^j x^k da_{jk} = 0  \notag \\
 \theta_i &= dp_{i} - \sum_j p_{ij} dx^j
   = da_i + \sum_j x^j da_{ij} = 0. \notag 
\end{align} 
Equivalently
\begin{equation}
\begin{pmatrix}
 2da_0 & da^t    \\
  da & dA   \\
\end{pmatrix} 
\begin{pmatrix}
1 \\
X\\
\end{pmatrix} = 0,   
\end{equation}
where $da^t = ( da_1, .. \, da_n )$, $( dA )_{ij} = da_{ij}
$, and $X^t = ( x^1, .. \, x^n )$, which is now considered
as a vector valued function on $\mathcal{S}$.
In other words, $\mathcal{S}$, as an $n$ 
-dimensional submanifold in the space of solutions
$\{ \,a_0, a_i, a_{ij} \, \}$,
is not only a null submanifold with respect to the 
symmetric $(n + 1)$ differential 
\[ \textnormal{det} 
 \begin{pmatrix}
 2da_0 & da^t    \\
  da & dA   \\
\end{pmatrix},
\]
but in fact it is a \emph{singular} null submanifold, 
meaning the matrix valued 1-form above has a null vector
as in (8).

Conversely, suppose $\mathcal{S}$ is an $n$ -dimensional 
singular null submanifold in the 
space of quadric hypersurfaces in $R^{n + 1}$.  
Generically, there is a vector valued function
$( 1, V^t )$  along $\mathcal{S}$ with 
$V^t = ( v^1, .. \, v^n )$ such that
\[ 
\begin{pmatrix}
  2da_0 & da^t    \\
  da & dA   \\
\end{pmatrix}
\begin{pmatrix}
 1\\
 V\\
\end{pmatrix}= 0,
\]
and $ dv^1 \wedge dv^2 \wedge .. \, dv^n \ne 0$
on $\mathcal{S}$. From the argument above, it is clear 
that the formula (7) with $x^i$ replaced by $v^i$ 
describes
a hypersurface in $R^{n + 1}$ that is the second order
developable to the given family of hyperquadrics $\mathcal{S}$.

Thus, at least in this \emph{flat} example, the vanishing
of the $(n+1)$ symmetric differential
\[ \textnormal{det}  
 \begin{pmatrix}
 2da_0 & da^t    \\
  da & dA   \\
\end{pmatrix}. 
\]
is a necessary condition for the contact of the neighboring 
Legendrian submanifolds. In case $n=1$, it is also
sufficient. We mention that for general
\emph{nonflat} family of Legendrian $n$-folds with
$n \geq 2$, 
the condition of contact of the neighboring submanifolds 
may not be as simple as this, as is discussed in [Ch2].

\newpage
 
\section{$G$ structure on the Space of Solutions}

\subsection{Contact of neighboring 
Legendrian Leaves}
 
The flat example considered above suggests the special 
class of the differential system $\mathcal{I}$ on $Z$ 
for which the conformal class of symmetric $(n+1)$ 
differentials, the vanishing of which represents a necessary
condition for the contact
of the neighboring Legendrian submanifolds, is well 
defined on the space of the leaves of the foliation $X$.
In fact, [Ch2] shows in case $n = 1$, the vanishing of a 
single relative invariant associated to $\mathcal{I}$ is 
both necessary and sufficient condition for a conformal 
class of quadratic  differential to be well defined 
on the space of solutions.

A higher dimensional analogue of this result exists and can
be described as follows. We continue to use the
notation adopted in Section \textbf{3}.

\begin{proposition}
Let $F_1 \to Z \to Y$ be the bundle associated to a
differential system $\mathcal{I}$ with a pseudo connection
such that the structure equation  (6) is true.
The conformal class of the symmetric $(n+1)$ differential
\begin{equation}
 \textnormal{det} \begin{pmatrix}
    2\theta_0 & \theta^t    \\
     \theta   & \Theta   \\
   \end{pmatrix}
\end{equation}
is well defined on the space of solutions if 
the bundle $F_1 \to Z$ admits a reduction to a subbundle
$F \subset F_1$ on which
\begin{align}
       T_0 &= 0,  \\
\pi^k_{ij} &\equiv \frac{1}{2} 
        ( \delta_{ik} \beta_j+\delta_{jk} \beta_i ) 
        \quad \mod \theta_0, \theta, \Theta, \notag \\
\epsilon &\equiv 0 \quad \quad \mod \theta_0, \theta, 
            \Theta.  \notag 
\end{align}
The structure equations on $F$ become
\begin{align}
d\theta_0 &= -\rho \wedge \theta_0 - 
             \theta^t \wedge \omega  \notag \\
d\theta &= - \beta \wedge \theta_0 -\alpha \wedge \theta 
             - \Theta \wedge \omega \notag \\
d\Theta &= -\frac{1}{2}( \beta \wedge \theta^t  
                       - \theta \wedge \beta^t )
           -( \alpha - \frac{\rho}{2} ) \wedge \Theta 
           + \Theta \wedge ( \alpha - \frac{\rho}{2} )^t
           + T \notag \\
d\omega &=  - \mu \wedge \theta_0 - \gamma \wedge \theta 
            +( \alpha^t - \rho ) \wedge \omega,
            \notag
\end{align}
where $T$ is quadratic in $\{ \theta_0, \theta, \Theta \}$
with $ T \equiv 0 \mod \theta_0, \theta$.
\end{proposition}

Before we begin the proof, we wish to give a interpretation
of the reduction procedure in local coordinates as in (5).
Once we get the structure equation (6) starting from the
representation (2) of $\mathcal{I}$, the torsion $T_0$ is
an expression involving $\{ x^i, u, p_i, p_{ij} \}$,
$\{ F_{ijk} \}$ and their derivatives, and the group 
variables $\{ \lambda,$ $A,$ $b,$ $e,$ $A_0,$ $A_i,$ $. .$
 $\;$ $\}$. 

First,
the reduction $T_0 = 0$ means we must be able to solve the
equation $T_0 = 0$ by expressing ${A_i}$ in terms of 
$\{ x^i, u, p_i, p_{ij} \}$, $\{ F_{ijk} \}$ and their 
derivatives, and $b$. Once we impose this relation back 
into the structure equation, we have 
\[ d\Theta = -\epsilon \wedge \theta_0 
           -\frac{1}{2}( \beta \wedge \theta^t  
                       - \theta \wedge \beta^t )
           -( \alpha - \frac{\rho}{2} ) \wedge \Theta 
           + \Theta \wedge ( \alpha - \frac{\rho}{2} )^t
           + T_{0}^{'} 
\]
where $T_0^{'}$ is a new torsion term. 

The second reduction
$\pi^k_{ij} \equiv \frac{1}{2} 
        ( \delta_{ik} \beta_j+\delta_{jk} \beta_i ) 
        \quad \mod \theta_0, \theta, \Theta $
means we must be able to express $e$ in terms of $\{ x^i, u, p_i, p_{ij} \}$,
$\{ F_{ijk} \}$ and their derivatives such that $T_0^{'}$
does not have any terms of the form $\theta_i \wedge 
\omega^j$. 
Finally, the third reduction 
$\epsilon \equiv 0 \, \mod \theta_0, \theta, \Theta \; $
means that when we impose the expression for $e$ obtained as
above back into the structure equations, the resulting
torsion $T$ should not have any terms involving 
$\omega^i$'s.
\vspace{1pc}

\emph{Proof of the Proposition}
Note that the structure equation on $F$ can be rewritten in a 
matrix form
\begin{align}
   d\begin{pmatrix}
     2\theta_0 & \theta^t \\
      \theta   & \Theta  
    \end{pmatrix} &= - 
    \begin{pmatrix}
     \frac{1}{2}\rho    & -\omega^t \\
     \frac{1}{2}\beta   & \alpha - \frac{1}{2}\rho 
    \end{pmatrix} \wedge
    \begin{pmatrix}
     2\theta_0 & \theta^t \\
      \theta   & \Theta  
    \end{pmatrix} \\
  &\quad +\begin{pmatrix}
     2\theta_0 & \theta^t \\
      \theta   & \Theta  
    \end{pmatrix} \wedge 
    \begin{pmatrix}
      \frac{1}{2}\rho  &  \frac{1}{2}\beta^t\\
     - \omega   & \alpha^t -  \frac{1}{2}\rho  
    \end{pmatrix} 
   + \begin{pmatrix}
      0 & 0 \\
      0 & T 
    \end{pmatrix}, \notag
\end{align} 
with
\[ d\omega =  - \mu \wedge \theta_0 - \gamma \wedge \theta 
            +( \alpha^t - \rho ) \wedge \omega.
\]
Let $G \subset Gl(N,R)$ be
the subgroup corresponding to the representation of 
$Gl(n+1,R)$ on the space of symmetric quadratic 
differentials on $R^{n+1}$ or on the space of 
$(n+1)$ by $(n+1)$ symmetric matrices.
The equation above, when $\theta_0, \theta,$ and  
$\Theta$ are interpreted as a basis of semibasic 1-forms,
is the equation of structure  
of an ordinary $G$ bundle over $X$. From the representation,
it is clear that the $G$ structure on $X$ 
induces a conformal class of symmetric $(n+1)$ differential 
(9) on $X$. $\square$

\vspace{1pc}
 
\begin{theorem}
Consider a Legendrian submanifold path geometry on $Z \to Y$ defined
by a Frobenius system $\mathcal{I}$ on $Z$. 
If the original $H_1 \subset Gl(N+n,R)$ bundle $F_1$ on
$Z$ admits a reduction to a subbundle $F$ on which the 
conditions (10) are satisfied, the bundle 
\[ F \to X \]
via $F \to Z \to X$ is an ordinary $G \subset Gl(N,R)$ 
structure on X with (11) as 
a structure equation. A generic $n$ -dimensional singular 
null submanifold in $X$ with respect to (9) corresponds to 
an $n$-parameter family of solutions of $\mathcal{I}$
that admits a developable Legendrian submanifold in $Y$.
\end{theorem}

Note that the notion of an $n$-plane in the tangent space
of $X$ being singular null with respect to the matrix valued
1-form  (9) is equivariant under the action of
$G$. It is likely the conditions (10) in the proposition
is also necessary for the symmetric differential to be
well defined on the space of solutions $X$.

\subsection{Dual Description}

In the present case, $X$ inherits a $G$ structure from being 
the space of Legendrian leaves. We briefly
discuss what could possibly be a description dual to 
this.  Let $X$ be an $N=\frac{1}{2}(n+1)(n+2)$ 
-dimensional manifold with a $G$ structure with the structure
equation (11) defined on the associated bundle $F_G \to X$.
Define an $m=\frac{1}{2}n(n+1)$-plane 
$\mathcal{N}_p \subset T_p X$ to be \emph{totally null} 
if it corresponds, under the identification of $T_p X$ with
the space of quadratic differentials on $R^{n+1}$ 
via the $G$ structure, to a subspace
\[  \mathcal{N}_p \cong \{ Q \in sym^2(R^{n+1*}) \, | 
        \,  v \, \lrcorner \, Q = 0 \; \mbox{for some nonzero}
    \;  v \in R^{n+1} \; \}. 
\]
Let $\mathcal{N} \to X$ be the
bundle of totally null $m$-planes with the associated 
projection $F_G \to \mathcal{N}$. The structure
group $G$ acts transitively on the set of totally null
$m$-planes, and from (11), we may take the differential 
system on $F_G$
generated by
\[ \{ \theta_0, \theta, d\theta \} 
\]
as the pullback of the differential system $\mathcal{J}_0$
on $\mathcal{N}$ whose solutions are totally null submanifolds. 

Consider a Frobenius differential system $\mathcal{J} 
\supset \mathcal{J}_0$ on $\mathcal{N}$ whose pull back 
to $F_G$ is generated by
\[ \{ \theta_0, \theta, \omega \}. 
\]
Equation (11) shows that the differential system above can indeed be
pushed down to $\mathcal{N}$. The definition of the 
bundle $\mathcal{N}$ and the second equation in (11)
also shows tha $\omega$, considered as a pseudo connection form
of the bundle $\mathcal{N} \to X$, measures the rate of 
change of the tangent $m$-planes along a totally null
submanifolds. The differential system $\mathcal{J}$
on $\mathcal{N}$ thus describes the canonical lifts 
to $\mathcal{N}$ of \emph{geodesic totally null} 
submanifolds in $X$. 

It is now clear from (11) that the space of leaves of 
the foliation of $\mathcal{N}$ by geodesic totally null 
submanifolds, $Y$, inherits a contact structure 
with an associated Legendrian submanifold path geometry. In fact  
equations (11) asserts that we can identify $\mathcal{N}$
with the bundle $Z$ of Legendrian $n$-planes over $Y$.

\begin{picture}(300,60)(-37,0)
\put(150,40){$Z \cong \mathcal{N}$} 
\put(180,25){$\searrow$}
\put(145,25){$\swarrow$}
\put(194,10){$X$}
\put(133,10){$Y$}
\put(192,27){$\pi_2$}
\put(133,27){$\pi_1$}
\end{picture}

Note the fibers of the bundle $Z \to Y$ project under $\pi_2$
to geodesic totally null submanifolds in $X$  and the fibers of the 
bundle $\mathcal{N} \to X$ in turn project under $\pi_1$ to
Legendrian submanifolds in $Y$ that give rise to the 
Legendrian submanifold path geometry.

We mention that for an arbitrary $G$ structure on an $N$ -dimensional 
manifold, the differential system $\mathcal{J}$ describing
geodesic totally null $m$-folds is in general not Frobenius. The 
Frobenius condition would force a single irreducible
piece of the torsion tensor associated to the $G$ 
structure to vanish.

\newpage

\section{Equivalence Problem}

In this section, we continue the analysis of the 
equivalence problem for the class of Legendrian submanifold path
geometry discussed in section \textbf{4}. The underlying idea that
guides us through the reduction procedure is the 
construction of an $\mathfrak{sl}(3,R)$ -valued projective connection
associated to the path geometry on a surface
[Ca1] and its generalization demonstrated by Chern in 
[Ch1] and [ChM].

We start with the structure equations (11) on $F$,
\begin{align}
d\theta_0 &= -\rho \wedge \theta_0 - 
             \theta^t \wedge \omega   \\
d\theta &= - \beta \wedge \theta_0 -\alpha \wedge \theta 
             - \Theta \wedge \omega \notag \\
d\Theta &= -\frac{1}{2}( \beta \wedge \theta^t  
                       - \theta \wedge \beta^t )
           -( \alpha - \frac{\rho}{2} ) \wedge \Theta 
           + \Theta \wedge ( \alpha - \frac{\rho}{2} )^t
           + T \notag \\
d\omega &=  - \mu \wedge \theta_0 - \gamma \wedge \theta 
            +( \alpha^t - \rho ) \wedge \omega,
            \notag
\end{align}
where $T=(T_{ij})$ is quadratic in $\{ \theta_0, \theta, \Theta \}$
with $ T \equiv 0 \mod \theta_0, \theta$. Explicitly,
we write 
\begin{align}
    T_{ij} &= T_{ji} \notag \\
        &=  \sum_k T_{ij}^k \theta_0 \wedge \theta_k 
     + \sum _{kl} T_{ij,kl} \theta_0 \wedge \Theta_{kl} 
 +  \sum_{kl} T_{ij}^{kl} \theta_k \wedge \theta_l 
 +  \sum_{klm} T_{ij,lm}^k \theta_k \wedge \Theta_{lm} \notag
\end{align}
with $\, T_{ij,kl} = T_{ij,lk}, \, T_{ij}^{kl} = - T_{ij}^{lk},$ 
and $\, T_{ij,lm}^k = T_{ij,ml}^k$.

Note that the pseudoconnection forms $\rho, \alpha,$ and 
$\beta$ in (12) are determined up to the change 
\begin{align}
\rho^* &= \rho + p \theta_0 \\
\alpha^{i*}_{j} &= \alpha^{i}_{j} + c^i_j \theta_0
              + \sum_k c^{i}_{jk} \theta_k \notag \\
\beta^{i*} &= \beta^i + c^i \theta_0 + \sum_k c^{i}_{k} \theta_k, \notag
\end{align}
where $p, c^i, c^i_j, c^{i}_{jk}=c^{i}_{kj}$ 
are independent variables. Following the procedure of the method of
equivalence [Ga], we shall determine the coefficients 
$c^i, c^i_j, c^{i}_{jk}=c^{i}_{kj}$ by imposing conditions
on the torsion $T$.

Applying the transformation (13) to the structure equation 
(12), we find
\begin{align}
   T_{ij}^{k*} &= T_{ij}^{k}
 -\frac{1}{2}( c^i \delta_{jk}+c^j \delta_{ik} ) \notag \\
  T_{ij}^{kl*} &= T_{ij}^{kl} 
-\frac{1}{2}(c^i_k \delta_{jl}- c^i_l \delta_{jk})
-\frac{1}{2}(c^j_k \delta_{il}- c^j_l \delta_{ik})\notag \\  
  T_{ij,kl}^* &= T_{ij,lk} 
-\frac{1}{2}(c^i_k \delta_{jl}+c^i_l \delta_{jk})
-\frac{1}{2}(c^j_k \delta_{il}+c^j_l \delta_{ik}) \notag \\
  & \quad + \frac{1}{2}p (\delta_{ik} \delta_{jl}+\delta_{il} \delta_{jk}) 
    \notag \\
T_{ij,lm}^{k*} &= T_{ij,lm}^k 
-\frac{1}{2}(c^i_{kl} \delta_{jm}+ c^i_{km} \delta_{jl})
-\frac{1}{2}(c^j_{kl} \delta_{im}+ c^j_{km} \delta_{il}). \notag
\end{align}
Upon contraction, it becomes (no summation convention)
\begin{align}
  T_{ii}^{i*} &= T_{ii}^{i} -c^i \notag \\
  T_{ii}^{ki*} &= T_{ii}^{ki} 
              - c^i_k, \quad i \ne k \notag \\  
  T_{ii,ii}^* &= T_{ii,ii}+p-2 c_{i}^i \notag \\
  T_{ii,im}^{k*} &= T_{ii,im}^k - c^i_{km}, 
            \quad i \ne m \notag \\
  T_{ii,ii}^{k*} &= T_{ii,ii}^k - 2 c^i_{ki}. \notag 
\end{align}
Hence, $c^i, c^i_j, c^{i}_{jk}=c^{i}_{kj}$ can be determined so as
to achieve (no summation convention)
\begin{align}
  T_{ii}^{i} &= 0  \\
  T_{ii}^{ki} &= 0 \notag \\  
  T_{ii,ii} &= 0 \notag \\
  T_{ii,im}^{k}+T_{ii,ik}^{m} &= 0,
  \quad i \ne m \; \; \mbox{and} \; \; i \ne k  \notag \\
  T_{ii,ii}^{k} &= 0, \quad i,m,k=1, \; . . ,\; n,\notag 
\end{align}
which we assume from now on.
The admissible transformations of the pseudo connection
forms preserving the equations (12) and the symmetry (14) of the
torsion $T$  now become
\begin{align}
\rho^* &= \rho + p \theta_0 \\
\alpha^{i*}_{j} &= \alpha^{i}_{j} + 
                   \delta_{ij} \frac{1}{2} p \theta_0\notag \\
\beta^{i*} &= \beta^i + \frac{1}{2} p  \theta_i \notag \\
\mu^{i*} &= \mu^{i} + h^i \theta_0 + \sum_k h_{ik} \theta_k
             -\frac{1}{2} p \omega^i \notag \\
\gamma^{*}_{ij} &= \gamma^{*}_{ji} \notag \\
         &= \gamma_{ij} + h_{ij} \theta_0. \notag     
\end{align}
with  new independent variables $h^i$ and $h_{ij}=h_{ji}$. 

Consider the bundle
\[ B_1 \to F \]
whose fiber at each point of $F$ is the space of pseudo
connection forms $\rho, \alpha, \beta, \mu, \gamma$ for
which the equations (12) and (14) are satisfied. Equations
(15) then give explicit formulas, in terms of the parameters
$p, h^i, h_{ij}=h_{ji}$, for the tautological forms
on $B_1$, which exist by the definition of the bundle 
$B_1$. We drop $^*$ and use $\rho, \alpha, \beta, \mu,
 \gamma$ to denote the corresponding tautological forms.
Set
\begin{align}
d\rho &=-\beta^t \wedge \omega - \theta^t \wedge \mu
         + \theta_0 \wedge \psi + \Omega_{\rho} \\
d\beta &=\rho \wedge \beta - \alpha \wedge \beta
  - \Theta \wedge \mu + \frac{1}{2}\theta \wedge \psi 
   + \Omega_{\beta} \notag \\
d\alpha &=-\alpha \wedge \alpha
-\frac{1}{2}(\beta^t \wedge \omega -\beta \wedge \omega^t)
-\frac{1}{2}(\theta^t \wedge \mu +\theta \wedge \mu^t) \notag \\
 & \quad \; \, \, - \Theta \wedge \gamma + \frac{1}{2}\theta_0 \wedge \psi
   + \Omega_{\alpha} \notag
\end{align}
where $\psi = -dp$. Exterior derivatives of the first
two equations in (12) then give structure equations for 
the tautological forms $\rho, \alpha, \beta$,
\begin{align}
&\Omega_{\rho} \wedge \theta_0 = 0 \\
&\Omega_{\beta} \wedge \theta_0 + \Omega_{\alpha} \wedge
 \theta + T \wedge \omega = 0. \notag
\end{align}
From the first equation in (17), we may arrange that
\begin{equation}
 d\rho =-\beta^t \wedge \omega - \theta^t \wedge \mu
         + \theta_0 \wedge \psi.  
\end{equation}
by modifying $\psi$ if necessary.
Note that at this stage, the equations (12), (14) and (18)
determine $\psi$ up to the change
\begin{equation}
\psi^* = \psi + t \, \theta_0 - \sum_k h^k \theta_k. 
\end{equation}
where $t$ is a new variable.
Put
\[ \Omega_{\beta}^i \equiv \sum_j P^i_j \theta_0 \wedge \theta_j
   + \sum_{jk} P^i_{jk} \theta_0 \wedge \Theta_{jk}
   + \sum_{jk} P^{i,jk} \theta_j \wedge \theta_{k} 
   + \sum_{klm} P^i_{k,lm} \theta_k \wedge \Theta_{lm} \]
mod $\omega, \rho, \alpha, \beta$, where $\Omega_{\beta}^t
= ( \Omega_{\beta}^1, .. \; \Omega_{\beta}^n) $, and 
$\, P^i_{jk}=P^i_{kj}$, $\,P^{i,jk}=-P^{i,kj}$, $\,
P^i_{k,lm}=P^i_{k,ml}$. 

Applying the transformation (15) and (19) to (16), we get 
\begin{align}
P^{i*}_j &= P^{i}_j - ( \frac{1}{4}p^2-\frac{1}{2}t)
              \delta_{ij} \notag \\
P^{i*}_{jk} &= P^{i}_{jk}
   +\frac{1}{2}(\delta_{ij}h^k+\delta_{ik}h^j) \notag \\
P^{i*}_{k,lm}&=P^{i}_{k,lm}
 -\frac{1}{2}(\delta_{im}h_{lk}+\delta_{il}h_{mk}). \notag
\end{align}
The contraction of the above gives (no summation convention)
\begin{align}
P^{i*}_i &= P^i_i -( \frac{1}{4}p^2-\frac{1}{2}t) \notag \\
P^{i*}_{ii} &= P^{i}_{ii}+h^i \notag \\
P^{i*}_{k,ii}+P^{k*}_{i,kk}&=
    P^{i}_{k,ii}+P^{k}_{i,kk} - 2h_{ik}. \notag
\end{align}
In terms of the bundle $B_1$, the above computations
imply that there exists a subbundle $B \subset B_1$ on which
(no summation convention)
\begin{align}
\sum_{i=1}^n P^{i}_i &= 0, \\
P^{i}_{ii} &= 0 \notag \\
P^{i}_{k,ii}+P^{k}_{i,kk}&=0, \quad i,k =1, \; .. \; n. \notag
\end{align}
In fact, (12), (14), (18) and (20) determine the pseudo connection forms 
$\rho, \alpha, \beta, \mu, \gamma$ and $\psi$ up to
the change 
\begin{align}
\rho^* &= \rho + p \theta_0 \notag \\
\alpha^{i*}_{j} &= \alpha^{i}_{j} + 
                   \delta_{ij} \frac{1}{2} p \theta_0\notag \\
\beta^{i*} &= \beta^i + \frac{1}{2} p  \theta_i, \notag \\
\mu^{i*} &= \mu^{i} - \frac{1}{2} p \omega^i \notag \\
\gamma^{*}_{ij} &= \gamma^{*}_{ji} \notag \\
         &= \gamma_{ij}  \notag   \\  
\psi^* &= \psi + \frac{1}{2} p^2 \theta_0. \notag
\end{align}
Note that $p$ is the fiber variable of the bundle 
$B \to F$. 

Now, the differential forms 
\[ \{ \theta_0, \theta, \Theta, \omega,
  \rho, \beta, \alpha, \psi, \mu, \gamma \} \]
are invariantly defined 
and form a basis of 1-forms on $B$. 
Set
\begin{align}
d\psi &= \rho \wedge \psi - ( \beta^t \wedge \mu
         + \mu^t \wedge \beta ) + \Omega_{\psi} \\
d\mu &=-\frac{1}{2}\psi \wedge \omega +\alpha^t \wedge \mu
         - \gamma \wedge \beta  + \Omega_{\mu} \notag \\
d\gamma &=-\rho \wedge \gamma+\frac{1}{2}( \mu \wedge 
  \omega^t - \omega \wedge \mu^t)+ (\alpha^t \wedge \gamma
        -\gamma \wedge \alpha ) + \Omega_{\gamma}. \notag
\end{align}
Then the structure equations so far can be written as
\begin{equation}
 d\Phi = -\Phi \wedge \Phi + \Omega 
\end{equation}
where $\Phi$ is the $\mathfrak{sp}(n+1,R) \subset \mathfrak{sl}(2n+2,R)$ -valued 
1-form
\begin{equation}
 \Phi = 
\begin{pmatrix}
\phi & \pi \\
\eta & -\phi^t \\
\end{pmatrix}
\end{equation}
with
\begin{equation}
 \eta = 
\begin{pmatrix}
2\theta_0 & \theta^t \\
\theta & \Theta \\
\end{pmatrix}, \; \;
\phi = 
\begin{pmatrix}
-\frac{1}{2}\rho  & -\frac{1}{2}\beta^t \\
\omega & -(\alpha^t-\frac{1}{2}\rho) \\
\end{pmatrix}, \; \;
\pi = 
\begin{pmatrix}
-\frac{1}{4}\psi & \frac{1}{2}\mu^t \\
\frac{1}{2}\mu & \gamma \\
\end{pmatrix},
\end{equation}
and the $\mathfrak{sp}(n+1,R) $ -valued curvature form
\begin{equation}
 \Omega = 
\begin{pmatrix}
\Omega_{\phi} & \Omega_{\pi} \\
\Omega_{\eta} & -\Omega^t_{\phi} \\
\end{pmatrix}
\end{equation}
with
\begin{equation}
 \Omega_{\eta} = 
\begin{pmatrix}
0 & 0 \\
0 & T\\
\end{pmatrix}, \; \;
\Omega_{\phi} = 
\begin{pmatrix}
0 & -\frac{1}{2}\Omega^t_{\beta} \\
0 & -\Omega^t_{\alpha} \\
\end{pmatrix}, \; \;
\Omega_{\pi} = 
\begin{pmatrix}
-\frac{1}{4}\Omega_{\psi} & \frac{1}{2}\Omega^t_{\mu} \\
\frac{1}{2}\Omega_{\mu} & \Omega_{\gamma} \\
\end{pmatrix}.
\end{equation}

Exterior differentiation of the second and the last
equation in (12) and (18) gives the following algebraic
equations satisfied by the curvature form.
\begin{align}
&\Omega_{\beta} \wedge \theta_0 + \Omega_{\alpha} \wedge
 \theta + T \wedge \omega = 0,  \\
&\Omega_{\mu} \wedge \theta_0 + \Omega_{\gamma} \wedge
 \theta + \Omega_{\alpha}^t \wedge \omega = 0, \notag \\
&\Omega_{\psi} \wedge \theta_0 - \Omega_{\mu}^t \wedge
 \theta + \Omega_{\beta}^t \wedge \omega = 0. \notag 
\end{align}
In particular, 
\begin{equation}
 \Omega \equiv 0 \; \mod \theta_0, \theta, \omega.
\end{equation}
Along with the structure equation (22), this implies that
each fiber of the bundle $B \to Y$ via $B \to F \to Z \to Y$
has the structure of the Lie group $P_1$, where 
$P_1 \subset Sp(n+1,R)$ is the stablizer of a line in 
$R^{2n+2}$ under the standard representation of $Sp(n+1,R)$.

We also mention that if $T=0$, the $G$-structure induced on $X$ is
torsion free, which by a result in [Br2] implies 
$\Omega = 0$. For this reason. we
call $T$ the \emph{primary invariant} of the Legendrian
path geometry under consideration.

\begin{theorem}
Given a Legendrian submanifold path geometry on  $Z \to Y^{2n+1}$ 
whose associated bundle $F_1$ admits a reduction to a
subbundle $F \to F_1$ with the structure equations $(12$,
there exists a $P_1 \subset Sp(n+1,R)$ bundle $B \to Y$ and a 
$\mathfrak{sp}(n+1,R)$
-valued 1-form $\Phi$ on $B$ given by $(23)$ with structure equations
$(22)$, or $(12)$, $(16)$, $(18)$ and $(21)$. The components of the curvature form
$\Omega$ in $(25)$ and $(26)$ satisfy $(27)$. At each point of $B$, $\Phi$
induces an isomorphism of the tangent space of $B$ with 
$\mathfrak{sp}(n+1,R)$. Two such Legendrian submanifold path geometries are equivalent
if and only if their associated bundles and the pseudoconnections 
are isomorphic.
\end{theorem}

A result of Cartan  [Ca2] implies the Legendrian submanifold path
geometries with the maximal dimension of symmetry vector 
fields are those for which the coefficients of the curvature
form $\Omega$ are all constants, the simplest being the case
when $\Omega=0$. The homogeneous Legendrian
path geometry realizing this \emph{flat} structure 
equation with the full group of symmetry $Sp(n+1,R)$ will
be examined in section \textbf{7}.

\newpage

\section{Normal Symplectic Connection}

It is well known in projective geometry that to every
torsion free affine connection on a manifold  there
exists a unique normal projective connection whose
paths coincide with the set of geodesics of the given affine
connection [Ch3].  Moreover, the sets of geodesics of two torsion 
free affine connections induce an equivalent path geometry
if and only if their associated normal projective 
connections are equivalent. The normal projective 
connection associated to a torsion free affine connection
thus captures the geometry of path defined by the set of
geodesics of the affine connection.

The purpose of this section is to generalize this idea 
to Legendrian submanifold path geometry and to discuss a 
\emph{Legendrdian connection} on a contact manifold that 
plays the role of a torsion free affine connection and induces
an associated Legendrian submanifold path geometry. Since  a \emph{path}
in Legendrian submanifold path geometry is an $n \geq 2$ -dimensional submanifold,
the integribility condition becomes nontrivial in this case.

We use $\theta_0$ to denote a (local) generator
of the contact structure on $Y$ or its pull back to the 
frame or other bundles over $Y$.

\subsection{Legendrian Connection}

Let $P \to Y^{2n+1}$ be the contact hyperplane vector bundle
of fiber dimension $2n$. 
\[ P \, = \, \{ \, ( \, p, \, v \, ) \; | \; p \in Y, \; 
  v \in T_p Y \; \mbox{with} \; \theta_0 ( v ) = 0 \, \} \]
Since
\begin{align}
 d(f \, \theta_0) &\equiv f \,  d \theta_0, \; \mod \theta_0,
                  \;   f \in C^{\infty}(Y) \notag \\
(d\theta_0)^n &\neq 0 \quad \mod \theta_0, \notag
\end{align}
the restriction of $d\theta_0$ induces a conformal 
symplectic structure on each fiber.
A basis $\{  A_i, \, B_i  \}$ of a contact hyperplane is 
called a conformal symplectic frame if
\begin{align}
d\theta_0 ( A_i, A_j ) &= 0,   \\
d\theta_0 ( B_i, B_j ) &= 0,  \notag \\
d\theta_0 ( A_i, B_j ) &= c \, \delta_{ij}, \; \; c \ne 0, 
              \quad i,j = 1,  \; .. , \, n.\notag 
\end{align}

Let $D$ be a conformal symplectic connection on the vector bundle $P$.
The infinitesimal displacement of a conformal symplectic frame field
is given by
\begin{equation}
 D\,( A, B ) = ( A, B ) \;\Psi, 
\end{equation}  
with $A=( A_1, \; . . , \; A_n )$, $B=( B_1, \; . . , \; B_n )$,
and $\Psi$ is a $\mathfrak{csp}(n,R)$ -valued connection form, i.e.,
\begin{equation}
\Psi= 
\begin{pmatrix}
\alpha + \rho I & \gamma \\
\Theta  & - \alpha^t + \rho I \\
\end{pmatrix} \notag
\end{equation}
where $\gamma, \; \Theta$ are symmetric $gl(n,R)$ -valued 
1-forms, $\alpha$ is a $gl(n,R)$ -valued 1-form, and 
$\rho$ is a scalar 1-form.
Under the change of the frame field 
\begin{equation} ( A^*, B^* ) = ( A, B ) \; g, 
\end{equation}
where $g$ is a $CSp(n,R)$ -valued function, we have
\[ D( A^*, B^* ) = ( A^*, B^* ) \; \Psi^*, \] 
with
\begin{equation}
\Psi^* = g^{-1} dg + g^{-1} \Psi \, g.
\end{equation}
The curvature form of the connection is defined by
\begin{equation}
 \Omega_{\Psi} = d \Psi + \Psi \wedge \Psi. \notag
\end{equation}
Differentiating (32), we have
\begin{equation}
 \Omega_{\Psi^*} = g^{-1} \Omega_{\Psi} g. 
\end{equation}

Take a local conformal symplectic frame field $\{  A_i, \, B_i \}$.
The identity of Cartan 
\[ d\theta_0(V_1, V_2) = V_1 \theta_0(V_2) - V_2\theta_0(V_1) 
- \theta_0([V_1,V_2]) \]
together with (29) gives
\begin{align}
\theta_0 ( [A_i, A_j] ) &= 0,   \\
\theta_0 ( [B_i, B_j] ) &= 0, \notag \\
\theta_0 ( [A_i, B_j] ) &= c \, \delta_{ij}, \; \; c \ne 0,  \quad i,j = 1, 
             \; .. , \, n.\notag 
\end{align}
Thus, if $V_1, \, V_2$ are vector fields tangent to the contact
hyperplane fields at each point spanning a isotropic 
plane field with respect to the conformal symplectic 
structure, then $[V_1,V_2]$ also lies in the contact hyperplane.
A connection $D$ is called \emph{isotropic torsion free} if
\[ D_{V_1} V_2 - D_{V_2} V_1 = [ V_1, V_2 ] \]
whenever the plane field generated by $V_1, V_2$ is isotropic
with respect to the conformal symplectic structure
on the contact hyperplane. We assume the connection $D$
to be isotropic torsion free from now on.

Let $F_P \to Y$ be the bundle of conformal symplectic
frames, which fits into the  diagram

\begin{picture}(200,72)(-37,-12)
\put(140,40){$F_P$} 
\put(159,30){$\searrow$}
\put(180,20){$Z$}
\put(159,11){$\swarrow$}
\put(144,20){$\downarrow$}
\put(142,-5){$Y$} 
\end{picture}

\noindent where $Z$ is the bundle of Legendrian $n$-planes.
The projection map $F_P \to Z$ is given by
\[ (\,  A_1, \, . . , \, A_n , B_1, \, . . , \, B_n \, )
   \to  A_1 \wedge A_2 \wedge \, . . , \, A_n \in Z \]
From the definition, there exists a set of tautological forms 
$\{ \omega^i, \theta_i \}$, $i=1, \, . . \, n$, conformal
symplectic coframe on $F_P$, defined
up to adding multiples of $\theta_0$. The equation
\begin{equation}
d\theta_0 \equiv - \sum_i \theta_i \wedge \omega^i
\quad \mod \theta_0
\end{equation}
in turn determines $\theta_0$ uniquely on $F_P$.

The isotropic torsion free 
condition on $D$ implies the following structure equation 
on $F_P$ dual to (30). 
\begin{equation}
d \begin{pmatrix}
    \omega \\
    \theta \\
\end{pmatrix} \equiv - \Psi \wedge  
\begin{pmatrix}
    \omega \\
    \theta \\
\end{pmatrix} + \sum_i^n t^i \, \theta_i \wedge \omega^i \; \;
\mod \theta_0
\end{equation}
where $\omega^t = (\, \omega^1, \, . . , \, \omega^n \, )$,
 $\theta^t = (\, \theta^1, \, . . , \, \theta^n \, )$, and 
$t^i$'s are functions on $F_P$. We call a isotropic
torsion free connection $D$ \emph{torsion free} if
we can modify $\omega$ and $\theta$ by adding 
multiples of $\theta_0$ to arrange
\begin{equation}
d \begin{pmatrix}
    \omega \\
    \theta \\
\end{pmatrix} = - \Psi \wedge  
\begin{pmatrix}
    \omega \\
    \theta \\
\end{pmatrix} + \theta_0 \wedge \Gamma
\end{equation}
with $\Gamma \equiv 0 \; \mod \; \omega, \, \theta$. 
Equation (37) determines the tautological
form  $\{ \omega, \; \theta \}$ uniquely on $F_P$.

Torsion free connections also admit the following equivalent but more
geometric description, which is a direct consequence
of (37).

\vspace{1pc}
\noindent \textbf{Definition} 
\textnormal{Torsion free connection } 

A conformal symplectic connection $D$ on the vector
bundle $P \to Y$ is called \emph{torsion free} if
there exists a line field in $Y$ transversal
to the contact hyperplane field such that for any vector
fields $V_1, \; V_2$ tangent to the contact hyperplane 
field,
\begin{equation}
 D_{V_1} V_2 - D_{V_2} V_1 = p([ V_1, V_2 ]) \notag
\end{equation}
where $p : TY \to P$ is the projection map induced
from the given line field. 

\vspace{1pc}
Of the torsion free connections, we consider the ones 
that give rise to a Legendrian submanifold path geometry on $Y$.
The natural analogue to the geodesics of the torsion free affine 
connection would be the geodesic Legendrian submanifolds, i.e., the
Legendrian submanifolds  whose tangent planes are
parallel along the submanifold with respect to the 
given connection. In terms of (37), they are the solutions
to the differential system $\mathcal{I}$ on $Z$ whose pull
back to $F_P$ is locally generated by
\begin{equation}
 \{  \theta_0, \theta, \Theta \}. 
\end{equation}
The structure equation (33) shows the differential system
(38) can indeed be pushed down to $Z$.

\vspace{1pc}
\noindent \textbf{Definition} 
\textnormal{Legendrian connection} 

A torsion free connection on the contact hyperplane
vector bundle $P \to Y$ is called \emph{Legendrian}
if the differential system $\mathcal{I}$ describing the
geodesic Legendrian submanifolds is
Frobenius on $Z$, or equivalently if the
differential system (38) is Frobenius on $F_P$. 

\vspace{1pc}
Set
\begin{align}
 \begin{pmatrix}
              \Omega_{\alpha+\rho } & \Omega_{\gamma} \\
              \Omega_{\Theta}  & - \Omega_{\alpha + \rho}^t 
               \end{pmatrix}  &= d\Psi + \Psi \wedge \Psi
             \notag \\
        &= \Omega_{\Psi}  \notag
\end{align}
Since the curvature form is quadratic in 
$\{ \theta_0, \, \theta_i, \, \omega^j  \}$,
the Frobenius condition corresponds to
\begin{equation}
 \Omega_{\Theta} \equiv 0 \quad \mod \theta_0, \, \theta.
\end{equation}
We assume the connection to be Legendrian from now on
and study the consequences of the equation (39). 

Let $V$ be the standard $2n$ -dimensional representation of
$Sp(n,R)$. Consider $S^2(V)$ $\otimes$ $\bigwedge^2(V)$
with the irreducible decomposition, [FH],
\begin{align}
 \mbox{$\bigwedge^2$}(V) &= \Gamma_{010..0} \oplus R  \\
 S^2(V) \otimes \Gamma_{010..0} &=
 \Gamma_{21} \oplus \Gamma_{20} \oplus \Gamma_{01} \quad 
 \mbox{if} \; n=2 \notag \\
&=\Gamma_{2100..0} \oplus \Gamma_{1010..0} \oplus \Gamma_{200..0} 
\oplus \Gamma_{010..0} \quad \mbox{if} \;n \geq 2. \notag
\end{align}
From (33) and (37), $\, \Omega  \mod \theta_0$ can be considered as a $Sp(n,R)$
equivariant $(S^2(V) \oplus R) \otimes \bigwedge^2(V)$ -valued
function on $F_P$. Equation (40)  shows that
in order for the equation (39) to be satisfied on $F_P$,
the curvature must be of the form
\begin{equation}
\Omega \equiv R \; \mbox{$\sum_{i=1}^n$} \theta_i \wedge 
          \omega^i \quad \mod \theta_0,
\end{equation}
where $R$ is a $\mathfrak{csp}(n,R)$ -valued function on $F_P$.

Exterior differentiation of (37) using (35) and (41) now gives
\[ \sum_i \theta_i \wedge \omega^i \wedge
   \Big( \Gamma + R\begin{pmatrix}
                  \omega \\
                  \theta \\
     \end{pmatrix} \Big) \equiv 0 \quad \mod \theta_0.\]
Since $n \geq 2$ and $\, \Gamma \equiv 0 \mod 
\theta, \omega$, we get
\begin{equation}
\Gamma = - R\begin{pmatrix}
                  \omega \\
                  \theta \\
              \end{pmatrix}. 
\end{equation}

Under the change (31), 
\[ \sum_i \theta_i \wedge \omega^i = \det(g)
   \, \sum_i \theta_i^* \wedge \omega^{i* }. \]
From (35) and (37), we get
\begin{equation}
d\theta_0 = -2 \rho \wedge \theta_0 - \sum_i \theta_i 
\wedge \omega^i + \theta_0 \wedge \sigma
\end{equation}
with $\sigma$ linear in  $\theta, \omega$. Differentiating
(43) and reducing $\mod \theta_0$ yields, 
\[ \sum_i \theta_i \wedge \omega^i \wedge \sigma
      \equiv 0 \quad  \mod \theta_0. \]
Thus $\sigma = 0 $  and we have
\begin{equation}
d\theta_0 = -2 \rho \wedge \theta_0 - \sum_i \theta_i 
\wedge \omega^i.
\end{equation}

By taking the exterior derivative of (44), 
\begin{equation} d \, \mbox{tr}(\Psi) = 2n \, d(\rho) 
  \equiv \mbox{tr}(R)\mbox{$\sum_i$} 
   \theta_i  \wedge \omega^i  \quad \mod \theta_0.
\end{equation}

\begin{proposition}
Let $Y^{2n+1}$ be a contact manifold with the contact line
bundle $L \subset T^* Y$, contact hyperplane bundle 
$P \to Y$, and its associated conformal symplectic frame 
bundle $F_P$. The equation
\textnormal{(35)} determines a section of the projection
$ L \oplus F_P  \to F_P $
and thus gives a well defined 1-form $\theta_0$
on $F_P$, which becomes a (local) generator of the
contact structure on $Y$ when pulled back by a section
of $F_P \to Y$. A Legendrian
connection on  $P$ 
induces a set of tautological forms
$\{ \, \theta, \omega \, \}$ with the structure equations 
\begin{align}
d\theta_0 &= -2 \rho \wedge \theta_0 - \mbox{$\sum_i$} 
\theta_i \wedge \omega^i \\
d \begin{pmatrix}
    \omega \\
    \theta \\
\end{pmatrix} &= - \Psi \wedge  
\begin{pmatrix}
    \omega \\
    \theta \\
\end{pmatrix} - \theta_0 \wedge R \begin{pmatrix}
    \omega \\
    \theta \\
\end{pmatrix} \notag \\
 d \Psi &\equiv - \Psi \wedge \Psi +
     R \; \mbox{$\sum_i$} \theta_i \wedge 
          \omega^i \quad \mod \theta_0. \notag
\end{align}
where $R$ is a $\mathfrak{csp}(n,R)$ -valued equivariant function on 
$F_P$.
A Legendrian submanifold path geometry is defined by the set of
geodesic Legendrian submanifold with respect to the 
given Legendrian connection, which, from the definition 
\textnormal{(40)}, is an $N = \frac{1}{2}(n+1)(n+2)$-parameter family
of Legendrian submanifolds in $Y$.
\end{proposition}
Note that (45) agrees with the trace of the  third equation 
in (46).

\subsection{Normal Symplectic Connection}

A projective connection on a manifold can be considered as
a family of torsion free affine connections all of which induce
 the same path geometry. In this section we wish to prove an analogous
result for Legendrian connections.

\begin{theorem}
Let $P \to Y^{2n+1}$ be the contact hyperplane bundle, and
$F_P$ the associated conformal symplectic frame bundle. Let
$D$ be a Legendrian connection on $P \to Y$. There exists
a bundle $B \to Y$ 

\begin{picture}(200,73)(-37,-12)
\put(140,40){$B$} 
\put(159,30){$\searrow$}
\put(180,20){$F_P$}
\put(159,11){$\swarrow$}
\put(144,20){$\downarrow$}
\put(142,-5){$Y$} 
\end{picture}

\noindent and an $\mathfrak{sp}(n+1,R)$ -valued 1-form $\Phi$, the normal
symplectic connection form, on $B$
with the following property.

1. The bundle $B \to Y$ has the structure of a principal
$P_1$ bundle, where $P_1 \subset Sp(n+1,R)$ is the stablizer of a 
line in $R^{2n+2}$ under the standard representation of
$Sp(n+1,R)$. 

2. At each point $u \in B$, $\Phi$ induces an isomorphism
between $T_u B$ and $\mathfrak{sp}(n+1,R)$. 

3. Under the right action by an element $g \in P_1$,
\[ \Phi_{u \, g} = Ad_{g^{-1}} \Phi_u. \]

4. A Legendrian connection $D'$ induces the same Legendrian
path geometry as $D$ if and only if $D'$ arises from a section
of the projection $B \to F_P$. 
\end{theorem}
The meaning of a Legendrian connection arising
from a section of the projection $B \to F_P$ will become
clear once the bundle $B$ is constructed.

Given a Legendrian connection $D$ with the structure equation
(46), we consider the space of Legendrian connections that
induce the same Legendrian submanifold path geometry equivalent 
as that of $D$. First of all, the deformation of the transversal
line field is expressed by
\begin{align}
 \omega^* &= \omega + p \, \theta_0 \\
 \theta^* &= \theta + q \, \theta_0 \notag
\end{align}
where $(p, q) = (p_1, . . , \, p_n, q_1, . . , \, p_n)$ is
a $R^{2n}$ -valued equivariant function on $F_P$. In order
to keep the equation (44), we must have
\begin{equation}
\rho^* = \rho + \frac{1}{2}( q^t \omega - p^t \theta )
         + s \, \theta_0
\end{equation}
where $s$ is a scalar equivariant function on $F_P$. 

Equation (47) and (48)  effect the second equation in (46). 
Equating (46) $\mod \theta_0$,  we have
\begin{equation}
\Psi.^* = \Psi. + \frac{1}{2} \begin{pmatrix}
p \\ q  \end{pmatrix} ( \theta^t, -\omega^t ) 
+ \frac{1}{2} \begin{pmatrix} \omega \\ 
\theta \end{pmatrix} ( q^t, -p^t )  +  x  + A \theta_0 
\end{equation}
where $\Psi = \Psi. + \rho I$, $A$ is a $\mathfrak{sp}(n,R)$ 
-valued equivariant function on $F_P$, and $x$ is in the 
kernel of the map $ \mathfrak{sp}(n, C)\otimes V^* \to V \otimes 
\bigwedge^2 V^*$. Here $V$ is the standard $2n$ -dimensional
representation of  $ \mathfrak{sp}(n, C)$.  
The kernel is isomorphic
to the irreducible representation $S^3(V)$
of $Sp(n,R)$  and,  since the ideal (38) describing the
Legendrian submanifold paths must be preserved, $x$ must be $0$.

Let $B_1 \to F_P$ be the bundle of all such point connections
with the explicit parametrization (47), (48), and (49).
We drop $^*$ and use the same notation to denote the tautological
forms on $B_1$. Taking the exterior derivative of 
(47), (48), and (49),
a computation shows that we can choose $A, \pi_0, \phi_0, \pi_0^0$ 
so as to achieve
\begin{align}
d\theta_0 &= -2 \rho \wedge \theta_0 - \mbox{$\sum_i$} 
\theta_i \wedge \omega^i  \\
d \begin{pmatrix}
    \omega \\
    \theta \\
\end{pmatrix} &= - ( \Psi. + \rho ) \wedge  
\begin{pmatrix}
    \omega \\
    \theta \\
\end{pmatrix} -  \theta_0 \wedge 
\begin{pmatrix}
    \pi_0 \\
   -\phi_0 \\
\end{pmatrix} \notag \\
d \Psi. &\equiv - \Psi. \wedge \Psi.
+ \frac{1}{2} 
\begin{pmatrix} 
\pi_0 \\ 
-\phi_0 \\
\end{pmatrix}
 \wedge ( \theta^t, -\omega^t ) + \frac{1}{2} 
\begin{pmatrix} 
\omega \\ 
\theta \\
\end{pmatrix}
 \wedge ( \phi_0^t, \; \pi_0^t ) \notag \\
 &\; \mod \theta_0  \notag \\
d \rho &= \frac{1}{2}( \theta^t \wedge \pi_0 + \omega^t 
\wedge \phi_0 ) + 2 \pi_0^0 \wedge \theta_0. \notag 
\end{align}
where
\begin{align}
 \pi_0 &\equiv dp  \notag \\
 \phi_0 &\equiv -dq  \notag \\
 \pi_0^0 &\equiv \frac{1}{2}ds \quad \mod \rho, \Psi., 
                   \theta_0, \theta, \omega. \notag
\end{align}

Let $B \subset B_1$ be the subbundle on which the set of
equations (50) hold. The 1-forms 
$\{ \theta_0, \theta, \omega, \Psi_., \rho \}$
are now uniquely defined on $B$, because of the fact that the only
solution $a \in gl(n,R)$ to the equation
\[ a \, \omega \wedge \theta^t + \omega \wedge \theta^t \, a = 0
\]
is $a=0$. The admissible change
of $\pi_0, \phi_0, \pi_0^0$ preserving (50) now becomes
\begin{align}
\pi_0^* &= \pi_0 + 2 y \, \theta_0 \notag \\
\phi_0^* &= \phi_0 + 2 z \, \theta_0 \notag \\
\pi_0^{0*} &= \pi_0^0 - \frac{1}{2}( \theta^t \, y + 
\omega^t \,  z ) + t \,  \theta_0 \notag
\end{align}
where $(y, z) \in R^{2n}$ and $t \in R$ are new variables.

Let us write
\[ d \Psi. = - \Psi. \wedge \Psi.
+ \frac{1}{2} \begin{pmatrix} \pi_0 \\ -\phi_0 \end{pmatrix}
 \wedge ( \theta^t, -\omega^t ) +
\frac{1}{2} \begin{pmatrix} \omega \\ \theta \end{pmatrix}
 \wedge ( \phi_0^t, \; \pi_0^t ) + \theta_0 \wedge \Upsilon 
\]
with $\Upsilon \equiv 0 \mod \theta, \omega$. Under the 
representation of $Sp(n, R)$,  $\Upsilon$ takes values in
$ S^2(V) \otimes V$ with the irreducible decomposition, [FH],
\[ \mbox{$S^2(V)$} \otimes V = \mbox{$S^3(V)$}
  \oplus V \oplus \Gamma_{110..0}. \]
Hence, we can absorb the $V$ piece by modifying 
$\pi_0, \phi_0$, and we have
\begin{equation}
\Upsilon \subset \mbox{$S^3(V)$} \oplus \Gamma_{110..0}
 \subset \mbox{$S^2(V)$} \otimes V. 
\end{equation}
This condition in turn determines $\phi_0$, and $\pi_0$ 
uniquely.

Put
\begin{align}
d\phi_0 &= (\rho + \alpha^t) \wedge \phi_0
          + 2\pi_0^0 \wedge \theta
          + \Theta \wedge \pi_0  + \Omega_{\phi_0}  \\
d\pi_0 &= (\rho - \alpha) \wedge \pi_0
          - 2\pi_0^0 \wedge \omega
          + \gamma \wedge \phi_0  + \Omega_{\pi_0} \notag 
\end{align}
In order to determine $\pi_0^0$, we write
\[ \Omega_{\phi_0}^i \equiv 
     \sum_j K^i_j \theta_0 \wedge \theta_j
 \mod \omega, \rho, \Theta, \alpha, \gamma,
\]
where 
$\Omega_{\pi_0}^t = ( \Omega_{\pi_0}^1, .. \; \Omega_{\pi_0}^n)$.
The 1-form $\pi_0^0$ is uniquely determined by imposing
\begin{equation}
 \sum_i K^i_i = 0. 
\end{equation}

Now, the structure equations (50), (51), (52) and (53) uniquely
determine $\{ \theta_0,$ $\theta,$ $\omega,$ $\Psi.,$ $\rho,$
 $\pi_0,$ $\phi_0,$ $\pi_0^0 \}$, and they form a basis
of 1-forms on $B$.  

Put
\[ d\pi_0^0 = 2\rho \wedge \pi_0^0 - \frac{1}{2}
        \phi_0^t \wedge \pi_0 + \Omega_{\pi_0^0}.
\]
Then the structure equations so far can be written as
\begin{equation}
 d\Phi = -\Phi \wedge \Phi + \Omega 
\end{equation}
where $\Phi$ is a $\mathfrak{sp}(n+1,R) \subset \mathfrak{sl}(2n+2,R)$ -valued 
1-form
\begin{equation}
 \Phi = 
\begin{pmatrix}
\phi & \pi \\
\eta & -\phi^t \\
\end{pmatrix}
\end{equation}
with
\begin{equation}
 \eta = 
\begin{pmatrix}
2\theta_0 & \theta^t \\
\theta & \Theta \\
\end{pmatrix}, \; \;
\phi = 
\begin{pmatrix}
-\rho  & -\frac{1}{2}\phi_0^t \\
\omega & \alpha \\
\end{pmatrix}, \; \;
\pi = 
\begin{pmatrix}
 \pi_0^0 & -\frac{1}{2}\pi_0^t \\
-\frac{1}{2}\pi_0 & \gamma \\
\end{pmatrix},
\end{equation}
and the curvature form
\begin{equation}
 \Omega = 
\begin{pmatrix}
\Omega_{\phi} & \Omega_{\pi} \\
\Omega_{\eta} & -\Omega^t_{\phi} \\
\end{pmatrix}
\end{equation}
with
\begin{equation}
 \Omega_{\eta} = 
\begin{pmatrix}
0 & 0 \\
0 & \Omega_{\Theta}\\
\end{pmatrix}, \; \;
\Omega_{\phi} = 
\begin{pmatrix}
0 & -\frac{1}{2}\Omega^t_{\phi_0} \\
0 & \Omega_{\alpha} \\
\end{pmatrix}, \; \;
\Omega_{\pi} = 
\begin{pmatrix}
 \Omega_{\pi_0^0} & -\frac{1}{2}\Omega^t_{\pi_0} \\
-\frac{1}{2}\Omega_{\pi_0} & \Omega_{\gamma} \\
\end{pmatrix}.
\end{equation}
Note that
\[ \begin{pmatrix}
 \Omega_{\alpha} & \Omega_{\gamma} \\
 \Omega_{\Theta} & -\Omega_{\alpha}^t
\end{pmatrix} = \theta_0 \wedge \Upsilon
\]
where $\Upsilon$ satisfies (51). 

It is easily verified 
$\Upsilon =0$ implies $\Omega = 0$. We call $\Upsilon$
the \emph{primary invariant} of the Legendrian submanifold path geometry
arising from a Legendrian connection. A curvature form
$\Omega$ for which the set of forms $\Upsilon$ 
and $\Omega_{\phi_0}$ satisfy (51) and (53) is called 
\emph{normal}.

The exterior differentiation of the last
equation in (50) gives the following algebraic
equations satisfied by the curvature forms.
\begin{equation}
\Omega_{\phi_0}^t \wedge \omega + \Omega_{\pi_0}^t \wedge
 \theta - 4 \Omega_{\pi_0^0} \wedge \theta_0 = 0. \notag 
\end{equation}
In particular, we have
\begin{equation}
 \Omega \equiv 0 \; \mod \theta_0, \theta, \omega, \notag
\end{equation}
which explains \emph{1.} in \textbf{Theorem 3}. The rest of the
statements in \textbf{Theorem 3} are also easily verified.

\subsection{Legendrian Flag}

From equation (46) or (50), consider the differential 
ideal
\begin{align}
\mathcal{I}_k &= \mathcal{I} \, \cup \, \{ \, \omega^{k + 1},
\, . . \, , \omega^n \, \} \, \cup \,
\{ \, \alpha^i_j \, | \, k+1 \leq i \leq n  \; \mbox{and}
 \; 1 \leq j \leq k \, \} \notag \\
 &\quad \mbox{for} \;  1 \leq k \leq n-1, \notag \\
\mathcal{I}_n &= \mathcal{I}. \notag
\end{align}
It is easy to check that each
 $\mathcal{I}_k$ is Frobenius on $B$. From (46),  $k$ 
-dimensional integral manifolds to $\mathcal{I}_k$ 
correspond to \emph{geodesic} 
isotropic $k$-folds in $Y^{2n+1}$ with respect to 
the given Legendrian connection, which from (50) are
invariantly defined independent of the choice of 
Legendrian connection of a normal symplectic 
connection. Since the structure group acts transitively 
on the set of isotropic $k$-planes, $\mathcal{I}_k$ defines
 a unique
geodesic isotropic $k$-fold tangent to each $k$-plane in a
 contact hyperplane  of $Y$
that is isotropic with respect to the induced conformal
symplectic structure.  

\begin{corollary}
Consider a Legendrian connection on a contact 
manifold $Y^{2n+1}$.
To each isotropic $k$-plane at a point of $Y$, there 
exists a unique geodesic istropic $k$-fold tangent to 
the given $k$-plane. The family of geodesic isotropic submanifold
paths defined by $\mathcal{I}_k$  
for $ 1 \leq k \leq n $ is in fact an invariant
of the associated normal symplectic connection.
\end{corollary}

Note that equation (54), when restricted to a geodesic
isotropic $k$-fold $\Sigma_k \subset Y$, gives rise to a bundle
\[ B_k \to \Sigma_k \]
with a $Gl(k+1, R)$ -valued Cartan connection form
$\phi_k$ given by the upper left hand corner $(k+1) \times (k+1)$ submatrix
of $\phi$. From (50) and (58), each fiber of this bundle has
the structure of a Lie subgroup $P_{k+1} 
\subset Gl(k+1, R)$, where
$P_{k+1}$ is isomorphic to the fiber of the bundle
\[Gl(k+1, R) \to RP^k.   \]

\newpage

\section{Examples}

\subsection{Flat Example}

Among the Legendrian submanifold path geometry is the simplest is that
of hyperquadrics in $R^{n+1}$, the \emph{flat} example
discussed earlier. It is easily verified this is the case
when the curvature form
\[ \Omega = 0 \]
in (22). The structure equation then suggests that the
Legendrian submanifold path geometry of hyperquadrics in $R^{n+1}$
is locally equivalent to the canonical homogeneous 
Legendrian submanifold path geometry on $RP^{2n+1}$.
 
Let $x^A, y^A, \; 0 \leq A \leq n$,  be the coordinates in 
$R^{2n+2}$ with the symplectic form
\begin{equation}
 \varpi = \sum_A dx^A \wedge dy^A. 
\end{equation}
Let $Q=Lag(n+1, R^{2n+2})$ be the space of Lagrangian 
$(n+1)$-planes in $R^{2n+2}$. The group 
of linear transformations $Sp(n+1,R) \subset Sl(2n+2,R)$
that preserves $\varpi$ then acts transitively on $Q$ and
the space of lines in $R^{2n+2}$, $RP^{2n+1}$, respectively.

In fact, it acts transitively via the product action 
on the incidence correspondence
\[ \mbox{I}= \{(l,E) \in RP^{2n+1}
           \times Q \; | \; l \subset E\}
\]
of dimension $(2n+1)+\frac{1}{2}n(n+1)$. The spaces
$\mbox{I}, \; RP^{2n+1}$ and $Q$, each of which is a homogeneous
space of $Sp(n+1,R)$, fit into the following diagram.

\begin{picture}(200,105)(-37,-4)
\put(120,80){$Sp(n+1,R)$} 
\put(150,62){$\downarrow$}
\put(150,44){I $=Sp(n+1,R)/P_1 \cap P_2$}
\put(135,26){$\swarrow$}
\put(160,26){$\searrow$}
\put(10,8){$Sp(n+1,R)/P_1=RP^{2n+1}$}
\put(170,8){$Q=Sp(n+1,R)/P_2$} 
\put(175,29){$\pi_2$}
\put(122,29){$\pi_1$}
\end{picture}

Here, we may choose $P_1 \subset Sp(n+1,R)$ to be the stablizer of the line
$l=\{(x^0, x^i=0, y^A=0)\}$ and $P_2 \subset Sp(n+1,R)$ to be the stablizer
of the Lagrangian $(n+1)$-plane $E=\{(x^A, y^A=0)\}$. 
Note that the fibers of the projections $\pi_1, \pi_2$
are $Lag(n, R^{2n})$ and $RP^n$ respectively.

The contact structure on $RP^{2n+1}$ can now be described as
follows. Take $l\in RP^{2n+1}$.  For a generator
$v\in R^{2n+2}$ of $l$, consider the 1-form 
\begin{equation}
 v \; \lrcorner \; \varpi. 
\end{equation} 
Since the generator $v$ is defined up to
a nonzero scalar multiple, the 1-form $v \; \lrcorner 
\; \varpi$ is also well defined on $T_l RP^{2n+1}$ up to 
a nonzero multiple. From (59), it is easy to see this is a 
contact structure as in (1). Also it follows from 
the construction that the fiber of the bundle 
$\mbox{I} \to RP^{2n+1}$ is the bundle of  Legendrian 
$n$-planes over $RP^{2n+1}$ with respect to the given contact
structure.

The description of the contact structure above naturally 
induces a Legendrian submanifold path
geometry on $RP^{2n+1}$. Take  $l\in RP^{2n+1}$ and
a Legendrian $n$-plane $E_0 \subset T_l RP^{2n+1}$. Then
there is a unique Lagrangian $(n+1)$-plane
$E \subset R^{2n+2}$ such that its image in $RP^{2n+1}$
under the projection $R^{2n+2}-\{0\} \to RP^{2n+1}$ is a $n$-
dimensional Legendrian submanifold of $RP^{2n+1}$ passing through
$l$ with $E_0$ as its tangent space. 

We wish to show that this path geometry is locally equivalent
to the geomerty of the hyperquadrics in $R^{n+1}$. Consider
a local coordinates parametrization of $RP^{2n+1}$ by
\[ \begin{pmatrix}
X^0 \\
X^i \\
Y^0 \\
Y^i 
\end{pmatrix}
= \begin{pmatrix}
1\\
x^i\\
2u-x^i p^i \\
p^i
\end{pmatrix}. \]
The contact structure is then locally generated by
\[ \sum_A X^A dY^A - Y^A dX^A = 2(du-\sum_i p^i dx^i). \]

Under these coordinates, the hyperquadrics 
\begin{align}
u &= a_0 + \sum_i a_i x^i + \sum_{ij} \frac{1}{2} 
                    a_{ij} x^i x^j  \notag \\
   p_i &= a_i + \sum_{j} a_{ij} x^j \notag 
\end{align}
with $a_{ij} = a_{ji}$ correspond to 
\begin{align}
Y^0 &= 2a_0 X^0 + \sum_i a_i X^i \notag \\
Y^i &= a_i X^0  + \sum_{j} a_{ij} X^j, \notag 
\end{align}
which are the equations in $R^{2n+2}$ that define  
Lagrangian $(n+1)$-planes.

\subsection{Example with Symmetry}

Let $M^n$ be a hypersurface in a space of constant 
curvature $\bar{M}^{n+1}_c$, $c = 1, 0,$ or $-1$, without
any continuous extrinsic symmetry, i.e., the subgroup of 
motions I$(M) \subset \mbox{Iso}(\bar{M}^{n+1}_c)$
that preserves $M$ is at most 
discrete. Let   $Y^{2n+1} = Gr(n, \bar{M}^{n+1}_c)  
\to \bar{M}^{n+1}_c$ be the bundle of hyperplane elements, 
which has a natural contact structure. Consider the
$N = \mbox{dim} \,\mbox{Iso}(\bar{M}^{n+1}_c) = \frac{1}{2}(n+1)(n+2)$-parameter
family of hypersurfaces in $\bar{M}^{n+1}_c$ generated by the 
action of $\mbox{Iso}(\bar{M}^{n+1}_c)$ on $M$. The cannonical lifts
of this family to $Y$ will in general be an $N$-parameter family 
of Legendrian submanifolds, nondegenerate 
in the sense discussed in Section \textbf{2}  if the second
fundamental form of $M$ is nondegenerate with distict and 
functionally independent eigenvalues. 

In order to see that this Legendrian submanifold path geometry is not flat,
it would suffice to show the nonexistence of injective
homomorphisms 
\begin{equation}
\mbox{Iso}(\bar{M}^{n+1}_c) \to Sp(n+1,R). \notag
\end{equation}
for $n \geq 1$.
 
\begin{lemma}
There does not exist an injective homomorphism
\begin{equation}
\textnormal{Iso}(\bar{M}^{n+1}_c) \to Sp(n+1,R). \notag
\end{equation}
for $n \geq 1$.
\end{lemma}

\emph{Proof.} $\,$ Let's take the case $c=1$,
\[
SO(n+1) \to Sp(n,R).
\]
with $n \geq 4$. The cases $n=2, 3$ and $c=-1, 0 \, $ follow 
from similar arguments. Also, we may assume the image 
of $SO(n+1)$ lies in the maximal compact subgroup 
$U(n) \subset Sp(n,R)$, in fact in $SU(n) \subset U(n)$.

The group $SO(n+1)$ has, as the first two representations of minimal 
dimensions, the standard representation $V$ of dimmension
 $(n+1)$, and the adjoint representation of dimension 
$\frac{1}{2}n(n+1) > 2n$ for $n \geq 4$. An injective
homomorphism $SO(n+1) \to SU(n)$ induces a faithful 
representation of $SO(n+1)$ of dimension $2n$, $W^{2n}$. 
From the injectivity and the inequality above, $W$ must
contain a $V^{n+1}$ piece. Since $SO(n+1) \subset SU(n)$,
it also preserves the complex conjugate $J(V)$, which
is impossible for dimension reasons. $\square$

\vspace{1pc}
The invariants of the Legendrian submanifold path geometry thus 
obtained in general are expressed  in terms of fourth 
order information on the orginal hypersurface.

\newpage

\section*{Bibliography}

\noindent 
[Br1] R. Bryant, 
  \emph{Two exotic holonomies in dimension four, path
geometries, and twistor theory},
   Proc. Symp. in Pure Math. \textbf{53} (1991), 33-88.

\vspace{1pc} 
\noindent 
[Br2] \underline{\quad}, \emph{Classical, exceptional,
and exotic holonomies: A status report}, Actes de la Table
de Geomtrie Differentielle en l'Honneur de Marcel Berger,
    Soc. Math. de France (1996), 93-166.

\vspace{1pc} 
\noindent 
[Br3] \underline{\quad}, \emph{Projectively flat Finsler 2-spheres
of constant curvature}, Selecta Math. 
(New Series 3) 3 (1997), 161-203.

\vspace{1pc} 
\noindent 
[Br4] \underline{\quad}, \emph{Elie Cartan and geometric
duality}, Journees Elie Cartan 1998 et 1999, Intitut Elie
   Cartan 16 (2000), 5-20. 
        
\vspace{1pc} 
\noindent 
[Br5] \underline{\quad}, \emph{Recent advances in the 
theory of holonomy}, (Expose 861), Seminaire Bourbaki,
   Volume 1998/99, Asterisque 266 (2000), 351-374.    
       
\vspace{1pc} 
\noindent 
[Ca1] E. Cartan, \emph{Sur les varietes a connexion
projective}, Quevres Completes, Partie 3, tome 1, 825-861.  

\vspace{1pc} 
\noindent 
[Ca2] \underline{\quad }, \emph{Sur la structure des groupes
infinis de transformations }, Qeuvres Completes, Partie 2,
 571-714.    

\vspace{1pc} 
\noindent 
[Ch1] S. S. Chern, 
   \emph{A generalization of the projective geometry
of linear spaces}, 
    S. S. Chern: Selected Papers \textbf{I}, 76-82.

\vspace{1pc} 
\noindent 
[Ch2] \underline{\quad }, 
   \emph{The geometry of the differential equation
  $y'''$ $=$ $F(x, y, y', y'')$}, 
    S. S. Chern: Selected Papers \textbf{I}, 45-59.

\vspace{1pc} 
\noindent 
[Ch3] \underline{\quad },
   \emph{On projective connections and projective 
relativity}, S. S. Chern: Selected Papers \textbf{III},
301-308.

\vspace{1pc} 
\noindent 
[ChM] S. S. Chern, J. K. Moser, \emph{Real hypersurfaces 
in complex manifolds}, Acta Math. 133 (1974), 219-271.

\vspace{1pc} 
\noindent 
[Ga] R. B. Gardner, \emph{The method of equivalence and 
its applications},  CBMS-NSF Regional Conference Series 
in Applied Mathematics, 58 (1989).

\vspace{1pc} 
\noindent 
[Fa] J. Faran, \emph{Segre families and real 
hypersurfaces}, Invent. Math. 60 (1980), no. 2, 135-172. 
 
\vspace{1pc} 
\noindent 
[FH] W. Fulton, J. Harris \emph{Representation theory. 
A first course}, Graduate Texts in Mathematics, 129. 
Readings in Mathematics. Springer-Verlag, New York (1991).

\vspace{1pc} 
\noindent 
[KN] S. Kobayashi, T. Nagano \emph{On projective 
connections}, J. Math. Mech. 13 (1964), 215-235.

\end{document}